\documentclass[letterpaper, 10 pt, conference]{ieeeconf}  % Comment this line out
\IEEEoverridecommandlockouts                              % This command is only
\pdfoutput=1

\usepackage{microtype}
\usepackage{slashbox}
\usepackage{graphicx}
\usepackage{subfigure}
\usepackage{booktabs} % for professional tables
\RequirePackage[loading]{tracefnt}
\usepackage{hyperref}
\usepackage{bbm}
\usepackage{cite}
\usepackage{dsfont}

\usepackage{amsmath}
\usepackage{amssymb}
\usepackage{mathtools}
\let\proof\relax

\usepackage{amsthm}

\usepackage{algorithm}
\usepackage{algpseudocode}
\newtheorem{thm}{Theorem}
\newtheorem{assumption}{Assumption}
\newtheorem{proposition}{Proposition}
\newtheorem{lemma}{Lemma}

\theoremstyle{definition}
\newtheorem{definition}{Definition}
\newtheorem{remark}{Remark}
\renewcommand{\qed}{\hfill\blacksquare}

\title{\LARGE \bf
Non-stationary Bandits with Habituation and Recovery Dynamics and Knapsack Constraints
}

\author{Qinyang He and Yonatan Mintz% <-this % stops a space
\thanks{The authors are with the Department of Industrial and Systems Engineering, University of Wisconsin-Madison, Madison, WI, 53705, USA {qhe57,ymintz}@wisc.edu}}

\begin{document}

\maketitle
\thispagestyle{empty}
\pagestyle{empty}

%%%%%%%%%%%%%%%%%%%%%%%%%%%%%%%%%%%%%%%%%%%%%%%%%%%%%%%%%%%%%%%%%%%%%%%%%%%%%%%%
\begin{abstract}

% Non-stationarity of the rewards and cost of each arm pulls are two important extension to the standard multi-armed bandit problem for modeling many real-world sequential decision making problems. One of the common non-stationary structures is the habituation and recovery effect and this effect is captured by the reducing or gaining unknown efficacy (ROGUE) bandit model. In this work, we study the combination of the ROGUE bandit and knapsack constraints so that the cost of arm pulls is also considered. We propose a pulling algorithm that determines the arm to pull by solving a linear program (LP) that takes the UCB estimates of the reward and cost as input . We show that the regret of our algorithm is sublinear in terms of the time horizon and total budget compared to the actual optimal pulling sequence. We also run experiments on a particular instance of ROGUE bandit with knapsacks to demonstrate the performance of our algorithm.
% 
% 
Multi-armed bandit models have proven to be useful in modeling many real world problems in the areas of control and sequential decision making with partial information. However, in many scenarios, such as those prevalent in healthcare and operations management, the decision maker's expected reward will decrease if an action is selected too frequently while it may recover if they abstain from selecting this action. This scenario is further complicated when choosing a particular action also expends a random amount of a limited resource where the distribution is also initially unknown to the decision maker. In this paper we study a class of models that address this setting that we call reducing or gaining unknown efficacy bandits with stochastic knapsack constraints (ROGUEwK). We propose a combination upper confidence bound (UCB) and lower confidence bound (LCB) approximation algorithm for optimizing this model. Our algorithm chooses which action to play at each time point by solving a linear program (LP) with the UCB for the average rewards and LCB for the average costs as inputs. We show that the regret of our algorithm is sub-linear as a function of time and total constraint budget when compared to a dynamic oracle. We validate the performance of our algorithm against existing state of the art non-stationary and knapsack bandit approaches in a simulation study and show that our methods are able to on average achieve a 13\% improvement in terms of total reward. 
\end{abstract}

%%%%%%%%%%%%%%%%%%%%%%%%%%%%%%%%%%%%%%%%%%%%%%%%%%%%%%%%%%%%%%%%%%%%%%%%%%%%%%%%
\section{INTRODUCTION}
Stochastic Multi-armed bandits (MAB) have become some of the most common models for analyzing problems in sequential decision making and control with partial information. MABs have been used to model various real life applications such as medical trials \cite{villar2015multi}, advertising \cite{schwartz2017customer}, and recommendation system \cite{zeng2016online}. In the classic stochastic MAB setting, a decision maker must select an action from a finite set of actions to maximize their long term reward without knowing \emph{a priori}  the reward distribution  associated with each action. This means that to be effective, their policy must balance choosing actions that may help them learn more about the system (exploration) with actions that given current information seem like they are likely to provide a high reward (exploitation). In general this reward distribution is assumed to be unchanging and stationary throughout the process and the decision maker is assumed to be able to take as many actions as they desire. However, in many real world scenarios such as those prevalent in personalized healthcare \cite{zhou2023spoiled} and operations management \cite{krishnasamy2021learning}, reward distributions may be non stationary and actions are restricted by resource constraints.
%
%
%Multi-armed bandits (MAB) are commonly used models for addressing sequential decision making problems. In the classic MAB setting, the decision maker chooses one action out of the possible action set at each time stop and receives a stochastic reward. The distribution for the reward provided by each action is initially unknown and the decision maker's goal is to maximize the total rewards in a given period. The universality of the MAB framework, capturing the exploration-exploitation trade-off, has led to its widespread applications, including medical trials \cite{c1}, advertising \cite{c2}, and recommendation system \cite{c3}. A diverse range of specific variants of the classic setting have been proposed and studied, with a fundamental variant being bandit learning in non-stationary environments, where reward distributions change over time. Several definition of non-stationarity has been extensively explored within the research community. \cite{c4}. 

Several frameworks have been used to analyze MABs with non-stationary rewards. The oldest model related to this problem is the restless bandit model proposed by Whittle \cite{whittle1988restless} where bandit rewards are dependent on internal transitioning states. In more recent literature there have been two  main families of non-stationary bandit models. The first family, are bandits where the non-stationary is not modeled explicitly but is subject to a total variation constraint \cite{garivier2008upper,yu2009piecewise,besbes2014stochastic,cheung2022hedging,liu2022non}. These approaches have been shown to be effective in settings where the reward distribution either changes very slowly over time or very infrequently. The other family of models considers structured non-stationarity \cite{mintz2020nonstationary, Gornet2022, Wei2018}. These approaches are more suitable for frequently changing rewards; however, they require the additional assumption that the decision maker has a model for how rewards can change over time. 

% In the seminal work by Whittle where the non-stationary bandits are referred as restless bandits, the rewards are dependent on arms' states which could change even if the arms are not pulled \cite{whittle1988restless}. Following that change-point or piecewise-stationary setting where the distribution of the rewards remains constant over epochs and changes at unknown
% time instants is studied \cite{garivier2008upper, yu2009piecewise}. In the more recent line of work, analyses are based on the variation budget, a quantification of how much the reward distribution could vary from different time steps \cite{besbes2014stochastic}, \cite{cheung2022hedging}, \cite{liu2022non}.

Several models in the literature have been devised for settings with bandit feedback and limited resources. \cite{tran2012knapsack} studies the case where the cost of pulling an arm is fixed and becomes known after the arm is pulled once and \cite{ding2013multi} extends the result to the scenario where the costs are random variables. A generalization of these two settings is known as Bandit with Knapsacks (BwK) problem where each arm pull consumes multiple constrained resources. In the standard BwK setting, \cite{badanidiyuru2018bandits} determines the policy by solving a sequence of linear programs (LP) that take the expected reward and resource consumption as input and the optimal values of these LPs are used for the regret analysis. 
%
%However, all of these models and methods are only devised for the case of stationary rewards.
Recently, models have been proposed to address the BwK case where rewards can be non-stationary \cite{liu2022non}. This model uses assumptions similar to those found in the literature that considers un-modeled non-stationarity that is bounded by a total variation constraint. This results in a dynamic regret bound that is $\mathcal{O}(\sqrt{T}\log T)$ with additional constants that depend on the total variation budget. Due to this assumption, this approach is well suited for cases where the distributions in the problem either change abruptly and infrequently or very slowly, and is not well suited for the case of frequently changing rewards. Thus, new approaches must be developed that can address both frequently changing non-stationary rewards and resource constraints. %both the reward distribution and consumption distribution may change over time but the change are bounded by an quantity called global non-stationarity budget. In the algorithm design, a evenly allocated budget over time is assumed and two LPs: one for step-wise arm pull decision, and the other for bounding optimal dynamic policy are used.

% Another natural extension to the standard MAB setting is assuming that pulling an arm is costly and the number of total possible arm pulls is limited by a fixed budget. \cite{tran2012knapsack} studies the case where the cost of pulling an arm is fixed and becomes known after the arm is pulled once and \cite{ding2013multi} extends the result to the scenario where the costs are random variables. A generalization of these two settings is known as Bandit with Knapsacks (BwK) problem where each arm pull incur consumption on multiply resources and each resource is budgeted. Under standard MAB with knapsack constraints setting, \cite{badanidiyuru2018bandits} uses a linear program (LP) that takes the expected reward and resource consumption to both drive the algorithm design and serve as the benchmark for regret analysis. \cite{liu2022non} further explores non-stationary bandit with knapsack constraints. In \cite{liu2022non} both the reward distribution and consumption distribution may change over time but the change are bounded by an quantity called global non-stationarity budget. In the algorithm design, a evenly allocated budget over time is assumed and two LPs: one for step-wise arm pull decision, and the other for bounding optimal dynamic policy are used.

In this paper, we propose methods  to address the challenges  of non-stationarity and knapsack constraints with frequently changing reward distributions. Our approach builds upon the literature related to structured non-stationarity and develops a new regret analysis for the case of dynamic regret. In particular we consider the case of non-stationarity where taking an action too frequently may reduce its reward, an effect known as habituation, while refraining from taking an action may increase its reward, known as recovery. Bandit models that address these effects are referred to as bandits with reducing or gaining unknown efficacy (ROGUE bandits) \cite{mintz2020nonstationary}. This form of non-stationarity is common in many real world applications such as personalized healthcare treatment \cite{liao2020personalized} and online advertising \cite{bertsimas2007learning}.  In this paper we present the model we call  the ROGUE bandits with knapsacks (ROGUEwK) problem and present an upper confidence bound (UCB) based approach for solving this problem we call the ROGUE knapsack UCB (ROGUEwK-UCB) algorithm. We show that our approach is able to achieve a dynamic regret bound of $\mathcal{O}(\sqrt{T}\log T)$ meeting the existing non-stationary and stationary BwK regret bound up to a log factor \cite{liu2022non,badanidiyuru2018bandits}. We conduct a computational experiment and show how our approach an outperform existing approaches by 13\% in terms of maximum reward.

% More specifically, we focus on an variant of restless bandit called reducing or gaining unknown efficacy (ROGUE) bandits proposed by Mintz \cite{mintz2020nonstationary}. An outstanding feature of the ROGUE bandit is that it captures habituation and recovery phenomena which is common in many real-world settings such as personalized healthcare interventions and targeted online advertising.  ROGUE bandits model this sort of non-stationarity effect by assuming state dependent rewards and change of states due to arm pulls. Combined with the knapsack constraints, this formulation can effectively model many sequential decision making scenarios where both the impact and cost of actions are taken into account.

This paper is organized as follows. In Section \ref{sec:prob_state}, we formulate the problem and discuss  relevant preliminary results. In Section \ref{sec:alg_descript}, we discuss two relevant optimization problems and describe the ROGUEwK-UCB algorithm. In Section \ref{sec:regret_analysis}, we provide a theoretical analysis for our algorithm. In Section \ref{sec:num_exper} we run experiments on a particular substantiation of the ROGUEwK bandit, namely the ROGUE logistic regression model with knapsacks, to demonstrate the performance of our algorithm. In Section \ref{sec:conclusion} we conclude the paper.

\section{Problem statement}
\label{sec:prob_state}
In the ROGUEwK problem, the decision-maker is given a fixed  finite set of arms $\mathcal{A}$ (with $|\mathcal{A}|=m$). At each round $t$, the decision maker must play one arm denoted by $a_{t}$. Each arm pull $a_t \in \mathcal{A}$ at time $t$ provides a stochastic reward $r_{a_t,t}$ that has a sub-Gaussian distribution $\mathcal{P}_{a_t,x_{a_t,t}}$ with expectation $\mathbb{E}[r_{a_t,t}]=g(x_{a_t,t})$ for a bounded function $g_{a_t}:\mathcal{X} \rightarrow \mathds{R}$. Each action $a\in \mathcal{A}$ has a state $x_{a,t}$ with nonlinear dynamics $x_{a,t+1} = h_{a}(x_{a,t},\pi_{a,t})$ where $ \pi_{i, t}=\mathds{1}\left[a_t=i\right], h:\mathcal{X} \times \mathbb{B} \rightarrow \mathcal{X}$ is a known dynamics function, and $\mathcal{X}$ is a compact convex set such that $x_{a,t} \in \mathcal{X} \ \forall a,t$ and $x_{a,0}$ is initially unknown for $a \in \mathcal{A}$. The maximum time horizon $T$ is finite and known in advance. We also assume that there are $d$ resources that each has a budget $B_j$. Without loss of generality, we assume $B_j=B$ for all $j$. Each arm pull $a_t\in \mathcal{A}$ at time $t$ incurs a stochastic consumption $c_{a_t,j,t}$ on resource $j$ with support in $[0,1]$ and denote the expected consumption matrix to be $\mathbf{\mathbf{C}}$ where $\mathbf{\mathbf{C}}_{ij}$ denotes the expected consumption on resource $j$ for arm $i$. The realizations of consumption $c_{a,j,t}$ are independently and identically distributed and in each round the resource consumption of the pulled arm is revealed to the decision maker. The interaction between the learner and the environment terminates at the earliest time $\tau$ when at least one constraint is violated, i.e. $\sum_{t=1}^\tau c_{a_t,j,t}>B$, or the time horizon $T$ is exceeded. The objective is to maximize the expected cumulative reward until time $\tau$, i.e. $\mathbb{E}\left[\sum_{t=1}^{\tau-1} r_{a_t,t}\right]$. We measure the performance of algorithm/policy $\Pi$ by its regret which is defined as:
$$
\operatorname{Reg}(\Pi, T):=\mathrm{OPT}(T)-\mathbb{E}\left[\sum_{t=1}^{\tau-1} r_{a_t, t} \mid \Pi\right].
$$
Here $\mathrm{OPT}(T)$ denotes the expected cumulative reward of the optimal dynamic policy given all the information on the initial state, reward and cost distributions.

\subsection{Technical Assumptions on ROGUEwK}
% To guarantee that the ROGUE bandits have provable properties, following assumptions are needed.
As part of our analysis we make the following technical assumptions. First is a set of assumptions introduced in \cite{mintz2020nonstationary} for the analysis of bandits with ROGUE non-stationarity.
\begin{assumption}
\label{assump1}
    $r_{a, t}$ are conditionally independent given $x_{a, 0}$ (or equivalently, the complete sequence of $x_{a, t}, a_t$).
\end{assumption}

This is a fairly mild assumption that is a non-stationary analogue to the classical MAB assumption of i.i.d rewards. Essentially it implies that at any two points in time $t, t^{\prime}$ such that $t \neq t^{\prime}$, $r_{a, t} \mid\left\{x_{a, t}, \right\}$ is independent of $r_{a, t^{\prime}} \mid\left\{x_{a, t^{\prime}}, \right\}$.
%
%This assumption states that for any two time points $t, t^{\prime}$ such that $t \neq t^{\prime}$ we have that $r_{a, t} \mid\left\{x_{a, t}, \right\}$ is independent of $r_{a, t^{\prime}} \mid\left\{x_{a, t^{\prime}}, \right\}$, and it is a mild assumption because it is the closest analogue to the assumption of independence of rewards in the stationary MAB.

\begin{assumption}
    For all $a \in \mathcal{A}$, the reward distribution $\mathcal{P}_{a,x}$ has a log-concave probability density function (p.d.f.) $p_a(r \mid x)$ for all $x \in \mathcal{X}$.
\end{assumption}
This assumption provides regularity for the reward distributions and many common distributions (e.g., Gaussian and Bernoulli) have this property.

Define $f(\cdot)$ to be $L$-Lipschitz continuous if $\left|f(x_1)-f(x_2)\right| \leq L\left\|x_1-x_2\right\|_2$ for all $x_1, x_2$ in the domain of $f$, and define $f$ to be locally $L$-Lipschitz on a compact set $\mathcal{S}$ if it has the Lipschitz property for all points of its domain on set $\mathcal{S}$. Our next assumption is on the stability of the above distributions with respect to various parameters.

\begin{assumption}
    The log-likelihood ratio $\ell(r ;  x^{\prime}, x)=$ $\log (p(r \mid  x^{\prime}) / p(r \mid x))$ of the distribution family $\mathcal{P}_{a,x}$ is locally $L_f$-Lipschitz with respect to $x$ on the compact set $\mathcal{X}$ for all values of  $x^{\prime} \in \mathcal{X}$, and $g$ is locally $L_g$-Lipschitz with respect to $x$ on the compact set $\mathcal{X}$.
\end{assumption}
This assumption guarantees that when two sets of parameters have similar values, the resulting distributions will be close to each other. We also introduce an additional assumption regarding the functional form of the reward distribution family.
\begin{assumption}
    The reward distribution $\mathcal{P}_{a,x_a}$ for all $x_a \in \mathcal{X}$ and $a \in \mathcal{A}$ is sub-Gaussian with parameter $\sigma$, and either $p(r \mid  x)$ has a finite support or $\ell(r ; x^{\prime}, x)$ is locally $L_p$-Lipschitz with respect to $r$.
\end{assumption}
This assumption is essential to guarantee that sample averages closely approximate their means, and it is upheld by various distributions (such as a Gaussian location family with a known variance). We impose the following conditions regarding the dynamics governing the state of each action.
\begin{assumption}
    The dynamic transition function $h$ is bijective and $L_h$-Lipschitz continuous such that $L_h < 1$.
\end{assumption}
This assumption ensures that there are no rapid changes in the states of each action, and it implies stability in the dynamics. The last assumption, drawn from most existing BwK literature, pertains to the total budget constraint of the problem.

\begin{assumption}[Linear Growth]
\label{assump:linear_growth}
The resource budget $B=b T$ for some $b>0$ .
\end{assumption}
This assumption is needed to quantify the relationship between $B$ and $T$ and this is a common assumption used in BwK literature \cite{liu2022non, immorlica2022adversarial}.

\subsection{Preliminaries: Concentration Inequalities}
Before diving into the details of the ROGUEwK-UCB, we introduce some preliminary results that we use in our analysis. UCB algorithms have been commonly used for stochastic non-stationary MAB \cite{auer2002finite,liu2022non,mintz2020nonstationary}. To construct the UCB and LCB we need to use appropriate concentration inequalities for the parameter estimates of the rewards and costs.
%
%In statioanry MAB, the key step of this class of algorithms is to construct a confidence bound for the reward of an action based on past observations using Chernoff-Hoeffding bounds \cite{wainwright2019high}. Comparably, for ROGUE bandits we need to construct a confidence bound for the parameters $x_{a,0}$ characterizing the distribution of each action $a \in \mathcal{A}$ which is harder than because the non-stationarity violates the standard assumptions of many concentration results. 
%
To analyze the concentration of the cost parameter estimates  we use the Azuma-Hoeffding Inequality in the following form:

\begin{lemma}[Azuma-Hoeffding's Inequality \cite{azuma1967weighted}]
\label{lemma:azuma}
Consider a random variable with distribution supported on $[0,1]$. Denote its expectation as $z$. Let $\bar{Z}$ be the average of $N$ independent samples from this distribution. Then, $\forall \delta>0$, the following inequality holds with probability at least $1-\delta$,
\begin{equation}
|\bar{Z}-z| \leq \sqrt{\frac{1}{2 N} \log (\frac{2}{\delta})} .
\end{equation}
 More generally, this result holds if $Z_1, \ldots, Z_N \in[0,1]$ are random variables, $\bar{Z}=\frac{1}{N} \sum_{n=1}^N Z_n$, and $z=\frac{1}{N} \sum_{n=1}^N \mathbb{E}\left[Z_n \mid Z_1, \ldots, Z_{n-1}\right]$.
\end{lemma}
For analyzing rewards, a limitation of this inequality is that it requires the random variables to be independent and in the context of learning, requires the use of unbiased estimators. These conditions are violated in the case of ROGUE rewards as due to the structure of the model maximum likelihood estimators (MLE), which are in general biased, will be more effective then unbiased ones. This necessitates the use of an alternative concentration inequality. The key of this concentration has to do with a quantity called the trajectory Kullbek-Liebler (KL) divergence that is defined as follows:

\begin{definition}[Definition 1 from \cite{mintz2020nonstationary}]
    For some input action sequence $\pi_1^T$ and arm $a \in \mathcal{A}$ with dynamics $h_a$, given starting parameter values $x_{a, 0} \in \mathcal{X}$, let $\mathcal{T}_a(T) \subset \{1,...,T\}$ be the set of times when action $a$ was chosen up to time $T$, then define the trajectory KL-divergence between these two trajectories with the same input sequence and different starting conditions as
\begin{equation}
\begin{aligned}
& D_{a, \pi_1^T}( x_{a, 0} \|  x_{a, 0}^{\prime}) =\sum_{t \in \mathcal{T}_a(T)} D_{K L}(\mathcal{P}_{a, x_{a, t}} \| \mathcal{P}_{ a,x_{a, t}^{\prime}}),\\
& \quad=\sum_{t \in \mathcal{T}_a(T)} D_{K L}(\mathcal{P}_{a, h_a^t(x_{a, 0})} \| \mathcal{P}_{a, h_a^t(x_{a, 0}^{\prime})}).
\end{aligned}
\end{equation}
\end{definition}
Where $h_a^k$ represents the functional composition of $h_a$ with itself $k$ times subject to the given input sequence, $\mathcal{P}_{a, x}$ is the probability law of the system under parameters $x$, and $D_{K L}$ is the standard KL divergence. Using this quantity we use the following concentration result:

\begin{thm}[Theorem 1 from \cite{mintz2020nonstationary}]
\label{theo:rogue_concen}
Let $x^*_{a,0}$ be the true initial state of an arm and let $\hat{x}_{a,0}$ be the MLE estimate for this parameter. That is,  $ \hat{x}_{a, 0}=\arg \min \{-\sum_{t \in \mathcal{T}_a(T)} \log p(r_t \mid  x_{a, t}): x_{a, t+1}
.=h_a(x_{a, t}, \pi_{a, t})\}$, where $\{r_t\}_{t \in \mathcal{T}_a}$ are the observed rewards for action $a \in \mathcal{A}$. Let $n_a(T) = |\mathcal{T}_a |$ denote the number of times arm $a$ is played up to time $T$. Then for $\alpha \in(0,1)$, with probability at least $1-\alpha$, we have
$$
\frac{1}{n_{a}(T)} D_{a, \pi_1^T}( x_{a, 0}^* \| \hat{x}_{a, 0}) \leq B(\alpha) \sqrt{\frac{\log (1 / \alpha)}{n_{a}(T)}},
$$
where $B(\alpha)=\frac{c_f(d_x)}{\sqrt{\log (1 / \alpha)}}+L_p \sigma \sqrt{2}$ and
$$\begin{aligned}
c_f(d_x)= 8 L_f \operatorname{diam}(\mathcal{X}) \sqrt{\pi}+48 \sqrt{2}(2)^{\frac{1}{a_x}} L_f \\
\cdot \operatorname{diam}(\mathcal{X}) \sqrt{\pi d_x}.
\end{aligned}$$
\end{thm}
Here the inclusions of the terms related to $B(\alpha)$ account for the bias in the MLE estimation.

\section{ROGUEwK-UCB Algorithm}
\label{sec:alg_descript}
In this section, we explain the details of the ROGUEwK-UCB algorithm. For the stationary BwK, the optimal dynamic policy can be computed by solving a LP that takes the mean reward and mean consumption vectors as input \cite{liu2022non,immorlica2022adversarial}. Also, the expected cumulative reward of this optimal dynamic policy is used for regret analysis. Similarly we introduce a nonlinear optimization that upper bounds the expected reward of the optimal dynamic policy in non-stationary case and use it to construct the ROGUwK-UCB algorithm.
\subsection{Relation of single step and multi-step problems}
 Let $\mathbf{x}_{t}=(x_{1,t}, x_{2,t},...,x_{m,t})$, $\mathbf{g}_t=(g_{1}(x_{1,t}),g_{2}(x_{2,t}),...,g_{1}(x_{m,1}))$. Define $NLP(\mathbf{x}_0,T, \mathbf{C})$ as
\begin{subequations}
\label{optimization:NLP}
    \begin{align}
    \max\limits_{\boldsymbol
    \pi_{t}} & \sum_{t=1}^{T}  \boldsymbol
    \pi^\top_{t} \mathbf{g}_{t},\\
\text{s.t. } &\\
&  x_{a,t} = h_a(\boldsymbol
    \pi_{a,t}, x_{a,t}), & \forall a \in \mathcal{A}, \forall t \in \{1,...,T\},\\
& \sum_{t=1}^{T} \mathbf{c}^\top_{j} \boldsymbol
    \pi_{t} \leq B, & \forall j \in \{1,...,d\},\\
& \boldsymbol
    \pi_{t} \in \Delta_m, & \forall t \in \{1,...,T\}.
    \end{align}
\end{subequations}
Here $\Delta_m \in \mathbb{R}^m$ is the $m$-dimensional unit simplex and $\mathbf{c}_{j}$ is the $j$th column of the matrix $\mathbf{C}$. Notice that the optimal value of $NLP(\mathbf{x}_0,T, \mathbf{C})$ is an upper bound on the expected cumulative reward of the optimal dynamic policy because it is the linear relaxation of the actual decision making problem which requires all the variables to be binary. This nonlinear optimization is hard to solve without specific structure in costs and state transition dynamics. A similar problem where the rewards have no state dependency has proven to be PSPACE-hard \cite{papadimitriou1994complexity}. We instead consider solving a LP to make step-wise decisions. Define the single-step optimization problem at time $t$ $LP(\mathbf{g}_t, \mathbf{C})$ with respect to $\{x_{a,t}\}_{a \in \mathcal{A}}$,  to be
\begin{subequations}
\label{optimization:SLP}
    \begin{align}
    \max\limits_{\boldsymbol
    \pi_{t}} &  \quad \boldsymbol
    \pi^\top_{t} \mathbf{g}_{t},\\
\text{s.t. } &\\
&\sum_{\mathcal{A}}  \mathbf{c}_{j}^\top \boldsymbol
    \pi_{t} \leq b, & \forall j \in \{1,...,d\},\\
& \boldsymbol\pi_{t} \in \Delta_m.
    \end{align}
\end{subequations}
The single-step LP problem can be interpreted as determining the optimal pulling distribution under a normalized resource budget $b$. The following proposition establishes the relationship between the global nonlinear optimization problem and step-wise linear program.
\begin{proposition}
\label{prop:step_nlp}
$NLP(\mathbf{x}_0, T, \mathbf{C}) \leq T \cdot LP(\mathbf{g}_0, \mathbf{C}) + L_g \frac{1}{1-L_h}\operatorname{diam}(\mathcal{X})$
\end{proposition}

\proof
Denote the optimal solution for $NLP(\mathbf{x_0}, T, \mathbf{C})$ as $\{\boldsymbol\pi^{\textit{NLP}}_t\}_{t=1}^{T}$ and the optimal solution for $LP(\mathbf{g}_0, \mathbf{C})$ as $\pi^*_{0}$.
\begin{subequations}
    \begin{align}
        &NLP(\mathbf{x}_0, T, \mathbf{C}) - T \cdot LP(\mathbf{x}_0, \mathbf{C}),\\
        &= \sum_{t=1}^{T} \sum_{a\in \mathcal{A}} {\pi_{a,t}^{\textit{NLP}}} g(h^t(x_{a,0})-\pi^{*}_{a,0} g(x_{a,0}), \\
        &\leq \sum_{t=1}^{T} \max\limits_{a \in \mathcal{A}} g(h^t(x_{a,0}))-g(x_{a,0}), \\
        &\leq \sum_{t=1}^{T}L_g \max\limits_{a \in \mathcal{A}}||h^{t}(x_{a,0})-x_{a,0}||_2, \\
        &\leq \sum_{t=1}^{T} L_g\max\limits_{a \in \mathcal{A}}||h^{t}(x_{a,0}) - h^{t}(x')||_2,\\
        &\leq \sum_{t=1}^{T} L_g L_h^t \operatorname{diam}(\mathcal{X})\leq  L_g \frac{1}{1-L_h}\operatorname{diam}(\mathcal{X}).
    \end{align}
\end{subequations}
Inequality 5d-5e comes from the fact that $h$ is a bijective function. $x_{a,0}=h(x')=h^2(x'')=...=h^t(x''')$ for some $x',x'',x''' \in \mathcal{X}$.
$\qed$

\subsection{Algorithm Details}
Next, we present the ROGUEwK-UCB algorithm as shown in Algorithm \ref{alg:ROGUE_KUCB}. We initialize the algorithm by pulling each arm once. After the first $m$ rounds, in every time step $t$, we first compute the MLE estimates of each arm's initial states and calculate the UCB for rewards based on the estimates of the states using Theorem \ref{theo:rogue_concen}. The idea behind the UCB is as follows: from Theorem \ref{theo:rogue_concen}, we know that with high probability, the true initial states are within a certain trajectory divergence from the MLE estimates, so we find the largest possible value of $g(x_{a,t})$ within the designated confidence radius. We also compute lower confidence bounds (LCB) for costs of each arm pull on different resources. Then we solve a single step LP problem which takes the UCBs and LCBs as input. The optimal solution to this LP is the probability distribution according to which the arm is going to be played and we pick an arm randomly following this distribution.

\begin{remark}
    While Algorithm \ref{alg:ROGUE_KUCB} is written such that it requires knowledge of time horizon $T$, if it is to be run indefinitely one could use the doubling trick \cite{besson2018doubling} and still preserve its statistical properties.
\end{remark}

% For computing the UCBs, we first calculate MLE of $x_{a,0}$ for each arm $a \in \mathcal{A}$ and set the UCB of arm $a$ to be be the maximum possible $g$ at time $t$ of another starting state which is within the allowed trajectory divergence from $\hat{x}_{a,0}$. The computation of the LCB is more standard which is the empirical mean minus the confidence width.
\begin{algorithm} 
\caption{ROGUEwK-UCB}
\label{alg:ROGUE_KUCB}
\begin{algorithmic}[1]
\Require Transition function $\{h_a\}$, reward function $\{g_a\}$
\For{$t\leq |\mathcal{A}|$}
\State Pick an arm $a$ that hasn't been chosen before
\EndFor
\For{$|\mathcal{A}| \leq t \leq T$}
\hspace{\algorithmicindent} \For{$a \in \mathcal{A}$}

\State \hspace{\algorithmicindent} Compute
\State$( \hat{x}_{a, 0})=\arg \min \left\{-\sum\limits_{s \in \mathcal{T}_a(t)} \log p(r_s \mid  x_{a, s})\right.$ \newline $\left. :x_{a, s+1}=h_a(x_{a, s}, \pi_{a, s}) \text { for } t \in\{0, \ldots, T\}\right\}$\newline
\State$g^{UCB}_{a,t} = \max_{ x_{a, 0} \in  \times \mathcal{X}}\left\{g( h_a^t(x_{a, 0})):\right.$ \newline
$\left.\frac{1}{n_{a}(t)} D_{a, \pi_a^t}(x_{a, 0} \|  \hat{x}_{a, 0}) \leq\right. \left.B(6mT^2) \sqrt{\frac{ \log (6mT^2)}{n_{a_s}(t)}}\right\}$ \newline
\State $c^{LCB}_{a,j,t} = \hat{c}_{a,j,n_a(t)} - \sqrt{{\frac{1}{2n_{a}(t)}} \log (12mdT^2)} \newline \quad \forall j \in \{1,...,d\}$  where $\hat{c}_{a,j,n_a(t)} = \frac{1}{n_a(t)}\sum_{s\in \mathcal{T}_a(t)}c_{j,s}$
\hspace{\algorithmicindent} \EndFor
\State Solve the single-step problem $LP(g^{UCB}_{t},\mathbf{C}^{LCB})$ and denote its optimal solution by $\boldsymbol\pi^*_{t} = (\pi^*_{1,t},...,\pi^*_{m,t})$
\State $\text { Pick arm } a_t \text { randomly according to } \boldsymbol\pi_t^* \text {, i.e., } \mathbb{P}(a_t=i)=\pi_{i,t}^*$
\State Observe reward $r_t$ and consumption $c_{j,t}, \forall j$.
\State Terminate if budget is exceeded
\EndFor
\end{algorithmic}
\end{algorithm}
\section{Regret Analysis for ROGUEwK-UCB}
\label{sec:regret_analysis}
In the following section, we present the analysis for the regret of the ROGUEwK-UCB algorithm. We first bound the difference between the maximum time horizon and the termination time when budget is exhausted.
\begin{proposition}
\label{prop:stop_time}
The following inequality holds with probability at least $1-\frac{1}{2T}$,
\begin{equation}
    T-\tau \leq \frac{1}{b}(1+4\sqrt{m}+\sqrt{\frac{1}{2}\log(12mdT^2)})(\sqrt{\frac{T}{2}\log(12mdT^2)}).
\end{equation}
\end{proposition}

\proof
From Hoeffding's inequality \cite{wainwright2019high} we have that for any $a \in \mathcal{A}$, $j \in {1,...,d}$, $t \leq \min \{\tau, T\}$, with probability at least $1-\frac{1}{6mdT^2}$, $|\hat{c}_{a,j,t} - \mathbf{C}_{a,j}|\leq \sqrt{\frac{1}{2n_{a_t}(t)}\log(12mdT^2)}$. Then
\begin{equation}
\begin{aligned}
&|\mathbf{C}_{a_t,j} - c^{LCB}_{a_t,j,t}| \\
= &|\mathbf{C}_{a_t,j} - (\hat{c}_{a_t,j,t}-\sqrt{\frac{1}{2n_{a_t}(t)}\log(12mdT^2)})|, \\
\leq &|\mathbf{C}_{a_t,j} - \hat{c}_{a_t,j,t}| + \sqrt{\frac{1}{2n_{a_t}(t)}\log(12mdT^2)}, \\
\leq &2\sqrt{\frac{1}{2n_{a_t}(t)}\log(12mdT^2)}. 
 \end{aligned}
\end{equation}
Then for all $t \leq \min{\{\tau, T}\}$ with probability at least $1-T\cdot m \cdot d\cdot \frac{1}{6mdT^2} = 1-\frac{1}{6T}$, 
\begin{equation}
\label{cost_eq1}
    \begin{aligned}
        &|\sum_{s=1}^{t}(\mathbf{C}_{a_s,j} - c^{LCB}_{a_s,j,s})| \leq \sum_{s=1}^{t}|\mathbf{C}_{a_s,j} - c^{LCB}_{a_s,j,s}|, \\
        &\leq \sum_{s=1}^{t} 2\sqrt{\frac{1}{2n_{a_s}(s)}\log(12mdT^2)}, \\
        &= \sum_{a \in \mathcal{A}} \sum_{s \in \mathcal{T}_{a}(t)}2\sqrt{\frac{1}{2n_{a_s}(s)}\log(12mdT^2)}, \\
        &= \sum_{a \in \mathcal{A}} \sum_{s=1}^{|\mathcal{T}_{a}(t)|}2\sqrt{\frac{1}{2s}\log(12mdT^2)},\\
        &\stackrel{(a)}{\leq} \sum_{a \in \mathcal{A}} 4\sqrt{\frac{|\mathcal{T}_{a}(t)|}{2}\log(12mdT^2)} \stackrel{(b)}{\leq}  4\sqrt{\frac{mT}{2}\log(12mdT^2)}
    \end{aligned}
\end{equation}
where $(a)$ comes from the fact that $\sum_{n=1}^N \frac{1}{\sqrt{n}} \leq 2 \sqrt{N}$ and $(b)$ follows from the Cauchy-Shwartz inequality. By Lemma \ref{lemma:azuma} with probability at least $1-\frac{1}{6T} \leq 1-\frac{1}{6mdT^2}$,
\begin{equation}
\label{cost_eq2}
    |\sum_{s=1}^{t}c_{j,s}-\mathbf{C}_{a_s,j}|\leq \sqrt{\frac{t}{2}\log(12mdT^2)} \leq \sqrt{\frac{T}{2}\log(12mdT^2)}.
\end{equation}
Denote $\mathbf{c}^{\textit{LCB}}_{j,t} = (c^{\textit{LCB}}_{1,j,t},...,c^{\textit{LCB}}_{m,j,t})$. Note that $\mathbb{E}[\mathbf{c}_{j,t}^{{LCB}^\top} \boldsymbol\pi^*_{t}] = \mathbb{E}[c^{LCB}_{j,a_t,t}]$ and $c_{a_t,j,t}^{LCB}\in [-\sqrt{\frac{1}{2}\log(12mdT^2)},1]$, then by Lemma \ref{lemma:azuma} we have with probability at least $1-\frac{1}{6T} \leq 1-\frac{1}{6mdT^2}$,
\begin{multline}
\label{cost_eq3}
%\begin{aligned}
    |\sum_{s=1}^{t}\mathbf{c}_{j,s}^{{LCB}^\top} \boldsymbol \pi^*_{s} - c^{LCB}_{a_s,j,s}|,\\
    \leq(1+\sqrt{\frac{1}{2}\log(12mdT^2)})(\sqrt{\frac{T}{2}\log(12mdT^2)}).
%\end{aligned}
\end{multline}
Combining  \eqref{cost_eq1},\eqref{cost_eq2},\eqref{cost_eq3}, with probability at least $1-3 \cdot \frac{1}{6T} = 1-\frac{1}{2T}$,
\begin{equation}
\label{cost_eq4}
\begin{aligned}
    &|\sum_{s=1}^{t}c_{j,s}-\mathbf{c}^{LCB^\top}_{j,s}\boldsymbol\pi^*_{s}| \leq |\sum_{s=1}^{t}c_{j,s}-\mathbf{C}_{a_s,j}|+, \\
    &|\sum_{s=1}^{t}(\mathbf{C}_{a_s,j} - c^{LCB}_{a_s,j,s})|  + |\sum_{s=1}^{t}c^{LCB}_{a_s,j,s} - \mathbf{c}_{j,s}^{{LCB}^\top} \boldsymbol\pi^*_{s}|,\\
    &\leq (1+4\sqrt{m}+\sqrt{\frac{1}{2}\log(12mdT^2)})(\sqrt{\frac{T}{2}\log(12mdT^2)}).
\end{aligned}
\end{equation}
Without loss of generality, we analyze the case when $\tau \leq T$. At termination time $\tau$, let $c_{j,t}$ denote the realized cost, then
\begin{equation}
    \sum_{t=1}^{\tau}c_{j,t} \geq bT,
\end{equation}
for some $j \leq d$.
From the fact that for all time $t$, $\boldsymbol \pi^*_{t}$ is a feasible solution to the problem $LP(g^{UCB}_{t},\mathbf{C}^{LCB})$, we have
\begin{equation}
    \sum_{t=1}^{\tau} {\mathbf{c}^{LCB}_{j,t}}^\top \boldsymbol \pi^*_{t} \leq b\tau,
\end{equation}
Combining this inequality with \eqref{cost_eq4}, we have with probability at least $1-\frac{1}{2T}$,
\begin{equation}
\begin{aligned}
\sum_{t=1}^{\tau}c_{j,t} \leq b\tau+(1+4\sqrt{m}+&\sqrt{\frac{1}{2}\log(12mdT^2)})\\&\cdot(\sqrt{\frac{T}{2}\log(12mdT^2)}).   
\end{aligned}
\end{equation}
Therefore we have
\begin{equation}
\begin{aligned}
    (1+4\sqrt{m}+\sqrt{\frac{1}{2}\log(12mdT^2)})(\sqrt{\frac{T}{2}\log(12mdT^2)})\\ \geq b(T-\tau),
\end{aligned}
\end{equation}
which yields the desired result. $\qed$

Next we bound the absolute difference between cumulative realized rewards and sum of optimal values of single step LPs.
\begin{proposition}
\label{prop:reward_prop}
    The following inequality holds for all $t \leq \min \{\tau, T\} $with probability at least $1-\frac{1}{2T}$.
\begin{equation}
\begin{aligned}
    &|\sum_{s=1}^{t}r_s-{g_s^{UCB}}^\top \boldsymbol \pi^*_s|\leq \sqrt{2T\sigma^2\log(12T)},\\
    &+(\frac{1}{1-L_h} + \sqrt{\frac{T}{2}\log(12T)})L_g\operatorname{diam}(\mathcal{X}).
\end{aligned}  
\end{equation}
\end{proposition}

\proof
    From the Hoeffding's inequality for sum of independent Sub-Gaussian random variables \cite{wainwright2019high}, we have that with probability at least $1-\frac{1}{6T}$,
    \begin{equation}
    \label{reward_bound_1}
    \begin{aligned}
        &|\frac{1}{t}\sum_{s=1}^{t}(r_s-g_{a_s})| \leq \sqrt{\frac{2}{t}\sigma^2\log(12T)}, \\
        \Longrightarrow &|\sum_{s=1}^{t}(r_s-g_{a_s})| \leq \sqrt{2t\sigma^2\log(12T)} \leq \sqrt{2T\sigma^2\log(12T)}.
    \end{aligned}
    \end{equation}
    Note that for any $s$, $\mathbf{E}[\boldsymbol
    \pi_{s}^{*^\top}\boldsymbol g_s^{UCB}] = \mathbf{E}[g_{a_s}^{UCB}]$, by Lemma \ref{lemma:azuma}, we have with probability $1-\frac{1}{6T}$,
    \begin{equation}
    \label{reward_bound_2}
        |\sum_{s=1}^{t}x_s^{*^\top}g_s^{UCB}-g_{a_s}^{UCB}| \leq L_g\operatorname{diam}(\mathcal{X})\sqrt{\frac{T}{2}\log(12T)}.
    \end{equation}
We denote
\begin{multline}
    x_{a_s,s}^{UCB} = \operatorname{arg}\max_{x_{a, 0} \in \mathcal{X}}\left\{g(h_a^t(x_{a, 0})):\right.\\
    \left.\frac{1}{n_{a}(t)} D_{a, \pi_1^t}(x_{a, 0} \| \hat{x}_{a, 0}) \leq\right. \left.B(6mT^2) \sqrt{\frac{ \log (6mT^2)}{n_{a_s}(t)}}\right\}
\end{multline}
By Theorem \ref{theo:rogue_concen} we have with probability at least $ 1-T \times m \times \frac{1}{6mT^2} =1-\frac{1}{6T}$ the following inequality holds
\begin{equation}
\label{reward_bound_3}
\begin{aligned}
&|\sum_{s=1}^{t}(g_{a_{s}}-g^{UCB}_{a_s})| \leq \sum_{s=1}^{t} |g(h^s(x_{a_s}))-g_{a_s}(h^s(x^{UCB}_{a_s,s}))|\\
&\leq \sum_{s=1}^{t}L_g|h^s(x_{a_s,s})-h^s(x^{UCB}_{a_s,s})|\\
&\leq \sum_{t=1}^{T} L_g L_h^t \operatorname{diam}(\mathcal{X}) \leq  L_g \frac{1}{1-L_h}\operatorname{diam}(\mathcal{X})
\end{aligned}   
\end{equation}
Then combine \eqref{reward_bound_1},\eqref{reward_bound_2}, \eqref{reward_bound_3} with probability at least $1-3 \times \frac{1}{6T} = 1-\frac{1}{2T}$,
\begin{equation}
    \begin{aligned}
        &|\sum_{s=1}^{t}r_s-{g_s^{UCB}}^\top \boldsymbol \pi^*_s| \leq |\sum_{s=1}^{t}(r_s-g_{a_s})| \\
        +&|\sum_{s=1}^{t}(g_{a_{s}}-g^{UCB}_{a_s})| + |g_{a_s}^{UCB}-\sum_{s=1}^{t}\boldsymbol \pi_s^{*^\top}g_s^{UCB}|\\
        \leq & \Big(\frac{1}{1-L_h} + \sqrt{\frac{T}{2}\log(12T)}\Big)L_g\operatorname{diam}(\mathcal{X}) \\ 
        + & \sqrt{2T\sigma^2\log(12T)}
    \end{aligned}
\end{equation}
$\qed$

Combining Propositions \ref{prop:step_nlp},\ref{prop:stop_time},\ref{prop:reward_prop} we can derive the final bound on regret.
\begin{thm}
Under Assumption \ref{assump1}–\ref{assump:linear_growth}, the regret of Algorithm \ref{alg:ROGUE_KUCB} is upper bounded as
\begin{equation}
\begin{aligned}
    &\operatorname{Reg}(\Pi^{\textit{ROGUEwK-UCB}}, T)\\
    \leq &\frac{1}{b}(1+4\sqrt{m}+\sqrt{\frac{1}{2}\log(12mdT^2)})(\sqrt{\frac{T}{2}\log(12mdT^2)})\\ 
    &\cdot LP(\mathbf{g}_0, \mathbf{C}) + (\frac{2}{1-L_h} + \sqrt{\frac{T}{2}\log(12T)})L_g\operatorname{diam}(\mathcal{X})\\
    +& \sqrt{2T\sigma^2\log(12T)} + LP(\mathbf{g}_0, \mathbf{C}) \\
    = &\mathcal{O}( \frac{1}{b}\sqrt{mT}\log(mdT)).
\end{aligned}   
\end{equation}     
\end{thm}
\proof
Let $\boldsymbol\pi^{LP}_0$ denote the optimal solution to $LP(\mathbf{g}_0, \mathbf{C})$
\begin{equation}
    \begin{aligned}
        &\sum_{t=1}^{\tau-1} {\boldsymbol g_{t}^{UCB}}^\top \boldsymbol \pi^*_t - (\tau-1)LP(\mathbf{g}_0, \mathbf{C}),\\
        =&\sum_{t=1}^{\tau-1} ({\boldsymbol g_{t}^{UCB^\top}} \boldsymbol \pi^*_t - \boldsymbol g^\top_0 \boldsymbol \pi^{LP}_0 )\\
        \leq & \sum_{t=1}^{\tau-1} \max\limits_{a \in \mathcal{A}} g(h^t(x_{a,0}))-g(x_{a,0}), \\
        \leq &\sum_{t=1}^{\tau-1}L_g \max\limits_{a \in \mathcal{A}}||h^{t}(x_{a,0})-x_{a,0}||_2,\\
        \leq &\sum_{t=1}^{\tau-1} L_g\max\limits_{a \in \mathcal{A}}||h^{t}(x_{a,0}) - h^{t}(x')||_2,\\
        \leq &\sum_{t=1}^{\tau-1} L_g L_h^t \operatorname{diam}(\mathcal{X}) \leq   \frac{L_g}{1-L_h}\operatorname{diam}(\mathcal{X}).
    \end{aligned}
\end{equation}
Then by Propositions \ref{prop:step_nlp},\ref{prop:stop_time},\ref{prop:reward_prop} we have probability at least $1-\frac{1}{T}$,
\begin{equation}
    \begin{aligned}
        &OPT-\sum_{t=1}^{\tau-1}r_t \leq NLP(\mathbf{x}_0, T, \mathbf{C}) - \sum_{t=1}^{\tau-1}r_t,\\
        =&(NLP(\mathbf{x}_0, T, \mathbf{C}) - \sum_{t=1}^{\tau-1} \boldsymbol g_{t}^{UCB^\top} \boldsymbol \pi^*_t) + \\ &(\sum_{t=1}^{\tau-1} {\boldsymbol g_{t}^{UCB^\top}} \boldsymbol \pi^*_t - \sum_{t=1}^{\tau-1}r_t),\\
        \leq &(T-\tau+1)LP(\mathbf{g}_0, \mathbf{C}) +L_g \frac{1}{1-L_h}\operatorname{diam}(\mathcal{X}), \\
        +&(\frac{1}{1-L_h} + \sqrt{\frac{T}{2}\log(12T)})L_g\operatorname{diam}(\mathcal{X}) + \\ & \sqrt{2T\sigma^2\log(12T)}, \\
        \leq &\frac{1}{b}(1+4\sqrt{m}+\sqrt{\frac{1}{2}\log(12mdT^2)})(\sqrt{\frac{T}{2}\log(12mdT^2)})\\
        & \cdot LP(\mathbf{g}_0, \mathbf{C})
        +(\frac{2}{1-L_h} + \sqrt{\frac{T}{2}\log(12T)})L_g\operatorname{diam}(\mathcal{X}) ,\\  + &\sqrt{2T\sigma^2\log(12T)} + LP(\mathbf{g}_0, \mathbf{C})=\mathcal{O}( \frac{1}{b}\sqrt{mT}\log(mdT)).
        \end{aligned}       
\end{equation}
Note that $\mathrm{OPT}(T)$ is of linear $T$ ($\mathrm{OPT}(T) \leq T\cdot L_g\operatorname{diam}(\mathcal{X})$) that transforms the high probability bound into the expectation bound. $\qed$

\begin{remark}
This result is significant because $\mathcal{O}(\frac{1}{b}\sqrt{mT}\log(mdT))$ is sublinear in $T$, and for fixed $T$, it is also sublinear in the total budget $B$ based on the relationship between $T$ and $B$ introduced in Assumption \ref{assump:linear_growth}. Regarding other state of the art results, our results matches the $\Omega(\sqrt{mOPT}\log(T))$ for stochastic BwK setting in \cite{badanidiyuru2018bandits} given that $OPT=\Theta(T)$ up to several log factors. However it is important to note that unlike the stationary setting, our result contains constants that are dependant on the non-stationarity of the system. Our result also matches \cite{liu2022non} where the non-stationarity of BwK is defined by the global non-stationarity budget.
\end{remark}

\section{Numerical Experiments}
\label{sec:num_exper}
In this section, we perform numerical experiments to demonstrate the effectiveness of ROGUEwK-UCB. We consider a dynamic generalized linear model (GLM) \cite{mccullagh2019generalized,filippi2010parametric} which can be interpreted as non-stationary generalizations of the classical (Bernoulli reward) stationary MAB \cite{gittins1979bandit,lai1985asymptotically,garivier2011kl}. The exact dynamics are as follows: the transition function $h(x_t, \pi_t) = A_{a}x_t + B_{a}\pi_t +K_{a}$ where $A_{a},B_{a},K_{a}$ are matrices/vectors of the correct size and and the rewards $r_{a,t}$ are Bernoulli with a logistic link function of the form $\mathbb{E}[r_{a,t}]=g_a(x_{a,t}) = \frac{1}{1+\text{exp}(-\alpha_{a}-\beta^\top_{a}x_{a,t})}$ where $\alpha$ is the vector of the correct size and $\beta_a \in \mathcal{R}$. We assume there are three resource constraints and the consumption distribution is uniform. We include three arms in this experiment whose dynamics and support of consumption distribution are shown in Table \ref{table:dynamics} and Table \ref{table:consumption}. Arm 1 has small habituation and recovery effects and is stable in its reward (indicated by high $k$ and $\beta$) but it has an unproportionally high cost for one of the resources. Arm 2 has moderate habituation and recovery effects and moderate consumption while Arm 3 has strong habituation and recovery effects and on average less consumption than Arm 2.

We compare our ROGUEwK-UCB algorithm against two other algorithms: the first is the naive UCB algorithm \cite{auer2002finite} that does not take non-stationarity and resource consumption into consideration; the second is the sliding window upper confidence bounds algorithm (SW-UCB) \cite{liu2022non} that can handle non-stationary in both rewards and costs by using a sliding window on the UCB estimates. For both algorithms we used the theoretically optiamlly derived hyper-parameters \cite{liu2022non}. We set the maximum time horizon $T$ to be 1,000 and test the cumulative rewards for budget from 10 to 300. Each of the candidate algorithms was replicated 10 times. 

\begin{table}[b]
\centering
\caption{Experimental Parameters for the dynamics of each Arm}
\label{table:dynamics}
\resizebox{0.8\columnwidth}{!}{\begin{tabular}{
lccccccc}
\hline Action & $x_0$ &  $A$ & $B$ & $K$ & $\alpha$ & $\beta$ \\
\hline 0 &  0.1 & 0.2 & -0.5 & 0.8 & 0.2 & 0.8 \\
       1 &  0.3 & 0.7 & -1.2 & 0.4 & 0.5 & 0.3 \\
        2 &  0.9 & 0.5 & -2.0 & 1.0 & 0.1 & 1.0 \\
\hline
\end{tabular}}
\end{table}

\begin{table}[b]
\centering
\caption{Support for the uniform consumption of each arm for each resource}
\label{table:consumption}
\begin{tabular}{l|lll} \hline
\backslashbox{Arm}{Resource} & 1 & 2 &3\\ \hline
0 & [0.1,0.2] & [0.6,0.8] & [0.3,0.5] \\ 
1 & [0.2,0.3] & [0.3,0.4] & [0.1,0.5] \\ 
2 & [0.2,0.3]&[0.2,0.4]&[0.1,0.3] \\ \hline
\end{tabular}
\end{table}

\begin{figure}
    \centering
    \includegraphics[width=0.85\columnwidth]{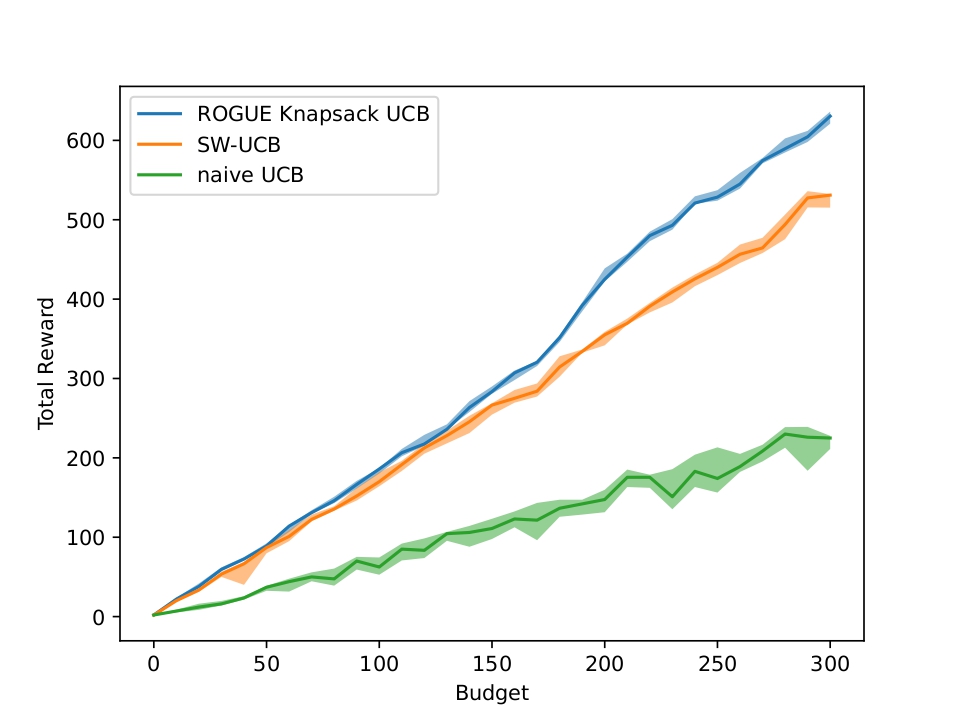}
    \caption{Cumulative Reward for each Algorithm}
    \label{fig:cumu_reward}
\end{figure}
\begin{figure}
    \centering
    \includegraphics[width=0.85\columnwidth]{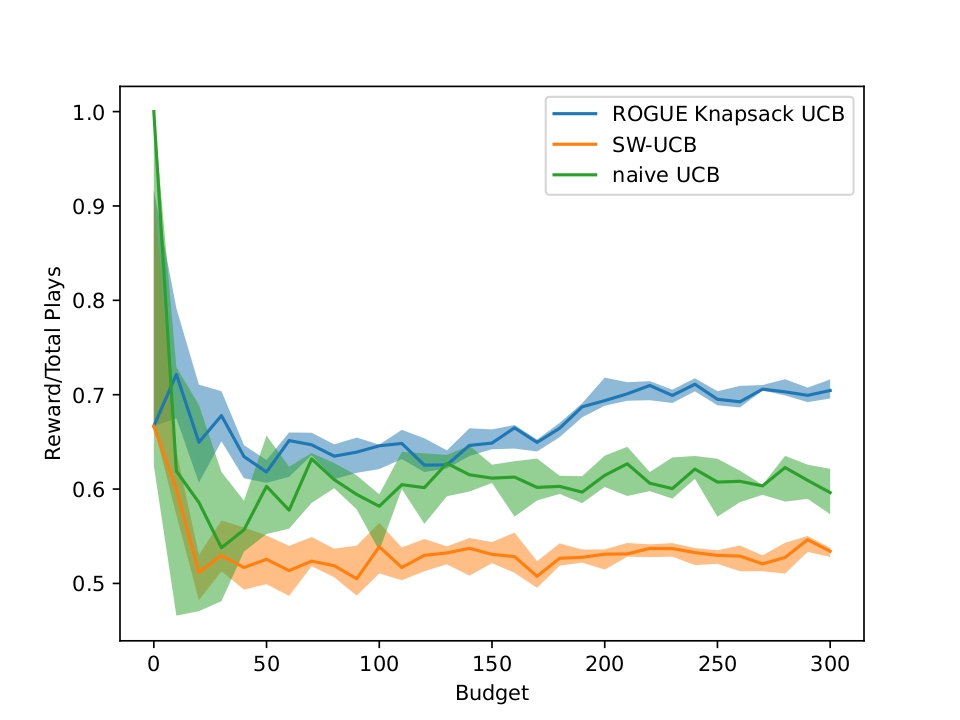}
    \caption{Average Reward per Play for each Algorithm}
    \label{fig:avg_reward}
\end{figure}
\begin{figure}
    \centering
    \includegraphics[width=0.85\columnwidth]{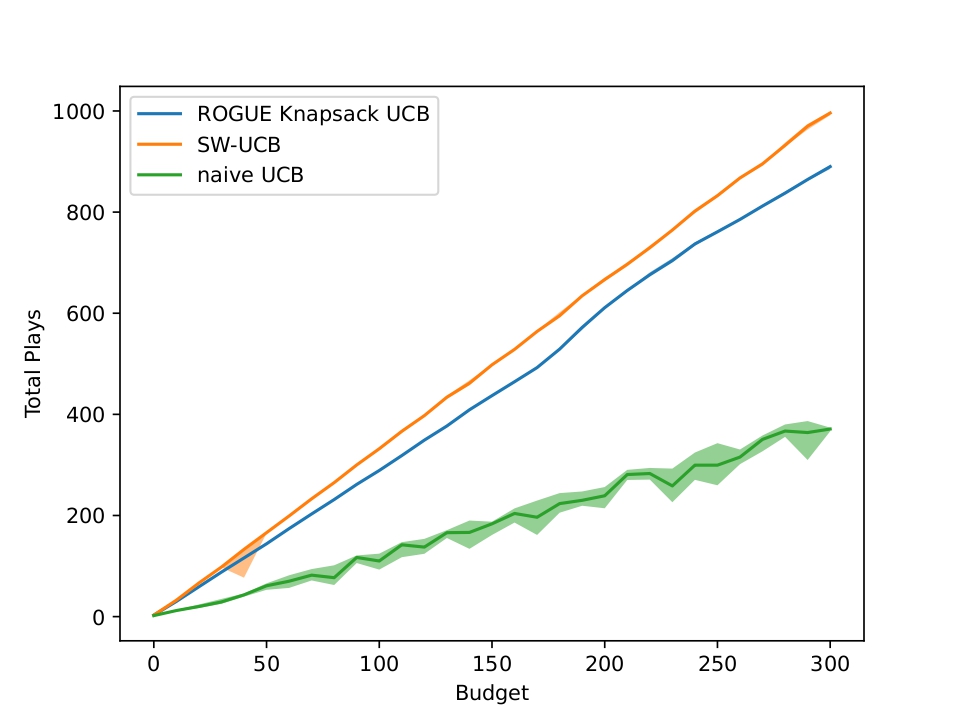}
    \caption{Total Plays for each Algorithm}
    \label{fig:plays}
\end{figure}
Figure \ref{fig:cumu_reward} shows the cumulative reward collected by all candidate algorithms within the maximum allowed time horizon; Figure \ref{fig:avg_reward} shows the average reward per play for each algorithm and Figure \ref{fig:plays} shows the total number of plays for each algorithms. The solid line represent the median value across 10 replicates and the shaded area represents the interquartile range of the values among 10 replicates. As depicted in the three plots, ROGUEwK-UCB achieves the most total reward across all the different budgets. In particular, compared with naive UCB, both BwK algorithms are cost-aware and avoid exhausting budget early by picking excessively costly arm and thus achieve higher cumulative reward. Compared with SW-UCB, although ROGUEwK-UCB has fewer plays before the budget is exhausted, it picks more cost-effective arms since it estimates reward based information of underlying non-stationary dynamics. ROGUEwK-UCB gains on average 13\% more total reward than SW-UCB across all different budgets.

\section{Conclusion}
\label{sec:conclusion}
In this work, we investigated non-stationary bandits with reducing or gaining unknown efficacy and knapsack constraints. We proposed an efficient UCB algorithm that determines the arms to play by solving a LP taking the UCB estimates of the rewards and LCB estimates of the costs as input. We showed that this algorithm achieves sublinear regret in terms of time horizon compared to the optimal dynamic oracle. Numerical experiments demonstrated that our algorithm outperforms other state-of-art algorithm for non-stationary bandits with knapsacks problems.

\bibliographystyle{IEEEtran}
\bibliography{IEEEabrv,main}

\end{document}